\documentclass[12pt,a4wide]{article}
\usepackage{amsmath, amsfonts, amssymb, amsthm, euscript, amscd, latexsym, mathrsfs, bm}
 \usepackage{titlesec}
\usepackage{amsmath, amsfonts, amssymb, amsthm, euscript, amscd, latexsym, mathrsfs, bm}
 \usepackage[pdftex,colorlinks=true,
                       pdfstartview=FitV,
                       linkcolor=blue,
                       citecolor=blue,
                       urlcolor=blue,
           ]{hyperref}

\def\aa#1{ \begin{align*} #1 \end{align*} }
\def\aaa#1{ \begin{align} #1 \end{align} }
\def\mm#1{ \begin{multline*} #1 \end{multline*} }
\def\mmm#1{ \begin{multline} #1 \end{multline} }

\newtheorem{thm}{\sc Theorem}
\newtheorem{lem}{\sc Lemma}

\newtheorem{rem}{\sc Remark}

\newcommand{\pl}{\partial}
\newcommand{\gt}{\geqslant}
\newcommand{\lt}{\leqslant}
\newcommand{\te}{\theta}
\newcommand{\sub}{\subset}

\newcommand{\dl}{\delta}
\newcommand{\al}{\alpha}
\newcommand{\gm}{\gamma}

 \newcommand{\sg}{\sigma}

\newcommand{\om}{\omega}
\newcommand{\mc}{\mathcal}
\newcommand{\Om}{\Omega}

\newcommand{\td}{\tilde}

\newcommand{\ox}{\otimes}

\renewcommand{\div}{\mathrm{div} \,}
\newcommand{\s}{{\alpha}}

\newcommand{\x}{\times}
\newcommand{\mto}{\mapsto}

\newcommand{\E}{\mathbb E}
\newcommand{\PP}{\mathbb P}

\newcommand{\C}{{\rm C}}

\newcommand{\ovl}{\overline}

\DeclareMathOperator{\ind}{\mathbb I}
\newcommand{\lap}{\Delta}
\newcommand{\nab}{\nabla}
\newcommand{\fdot}{\,\cdot\,}
\def\Rnu{{\mathbb R}}

\def\Qnu{{\mathbb Q}}
\def\Tor{{\mathbb T}^n}

\def\com#1{}

\long\def\symbolfootnote[#1]#2{\begingroup%
\def\thefootnote{\fnsymbol{footnote}}\footnote[#1]{#2}\endgroup}
\titleformat{\section}[hang]{\large\bfseries}{\thesection.}{1ex}{}{}
\titleformat{\subsection}[hang]{\normalsize\bfseries}{\thesubsection}{2ex}{}{}
\titleformat{\subsubsection}[hang]{\small\bfseries}{\thesubsubsection}{2ex}{}{}

\include{srctex.sty}

\begin{document}

 \author{Ana Bela Cruzeiro$^{1,2}$ and Evelina Shamarova$^3$
 }

 \date{}

 \title{ \Large
  \textbf{
 On a forward-backward stochastic system associated to the Burgers equation 
  }
 }

 \maketitle

 \vspace{-11mm}

 {\small
  \begin{center}
  \begin{tabular}{l}
  { $^1$\! Dep. de Matem\'{a}tica, IST-UTL.}  \\
  { $^2$\! Grupo de F\'isica Matem\'atica da Universidade de Lisboa.} \\
  { \hspace{2mm} E-mail: 
  \href{mailto:abcruz@math.ist.utl.pt}{abcruz@math.ist.utl.pt} }\\ 
  { $^3$\! Centro de Matem\'atica da Universidade do Porto.} \\
  { \hspace{2mm} E-mail: 
  \href{mailto:evelinas@fc.up.pt}{evelinas@fc.up.pt}}
  \end{tabular}
 \end{center}
 }

 \vspace{3mm}

 \begin{abstract}
We describe a probabilistic construction of  $H^{\alpha}$-regular solutions for the spatially periodic forced Burgers equation
by using a characterization of this solution through a forward-backward stochastic system.
 \end{abstract}
\section{Introduction}
Burgers equation, given by
$$
\pl_s y(s,\te) + (y,\nab)y(s,\te) - \nu \lap y(s,\te) + F(s,\te) = 0
$$
is sometimes presented as a simplified model for turbulence and describes the motion of
a compressible fluid with viscosity $\nu$  under the influence of a force $F$.
In this paper we establish a connection between a time-changed spatially periodic Burgers equation and a 
forward-backward stochastic system on the group of diffeomorphisms of a torus, 
similar to the characterization we have studied in [C-S] for the incompressible (Navier-Stokes equation) case.

It is well known that forward-backward systems are closely related to partial differential equations. 
For references in the subject one can use  for example those in  [D] or in the 
more recent work [C-S-T-V]. One difference in our approach is that the stochastic processes are defined  
in the group of diffeomorphisms of the underlying configuration space of the fluid (in our case, the torus) 
and not in the configuration space itself. This group is endowed with a Sobolev topology.
The importance of working with infinite dimensional geometry, in the line of thought
introduced by V. Arnold ([A]) for the Euler equation is, partly, that it allows 
to construct solutions which are ``automatically" Sobolev regular in the space variable. 
Also, which is more important, the geometric objects defined in the (infinite dimensional) 
path spaces may allow to prove several properties of the dynamics, 
such as stability for the corresponding flows.

Generally speaking, our approach can be regarded as  a stochastic deformation of the Lagrangian picture
keeping the ``mean velocity" (the Eulerian picture) of the motion unchanged. 
This means that the mean velocity, given by the drift, is still the relevant 
deterministic solution of the (velocity) equations of motion. 
This approach finds its roots in the works [N-Y-Z], [Y]. It is the point of view described in
[G], but is completely different from the
approaches that consist in perturbing the Eulerian (velocity) dynamics by a random noise.

In [C-S] we have developed in this spirit a  construction for the Navier-Stokes equations. 
We  derived a solution of the stochastic system from a solution of Navier-Stokes 
equation and, ``reciprocally", defined  Navier--Stokes solutions from the stochastic forward-backward flows. 
The incompressibility condition there makes the geometry much more delicate to study then in the present Burgers case and, in this respect, the 
Burgers equation is a simplification of the framework of [C-S]. On the other hand here we  
prove  an existence result of the stochastic forward-backward equation (without assuming the existence of the p.d.e. solution), 
a result which is not proved for the Navier--Stokes case. This is therefore the main result of this paper and the one that really distinguishes it from [C-S]. The methods we use to prove this result are close to those of Delarue in [D]. 
Therefore we obtain here  a completely probabilistic construction for the Burgers solutions.
We treat the torus case  since it is one of the simplest compact manifolds; the results, with the necessary adjustments, should extend to other manifolds.

We refer to [C-S] and references therein for the  general framework
of the stochastic approach to  partial differential equations such as Burgers and Navier--Stokes that we are dealing with.

\section{The main result}

Let us consider the spatially periodic backward Burgers equation in $\Rnu^n$:
\aaa{
\label{brgs}
\begin{cases}
\pl_s y(s,\te) + (y,\nab)y(s,\te) + \nu \lap y(s,\te) + F(s,\te) = 0\\
y(T,\te) = h(\te).
\end{cases}
}
It is obtained from the classical Burgers equation by means of the substitution
$u(t,\te) \leftrightarrow -u(T-t,\te)$, where $T>0$ is fixed arbitrary.
We assume that $h$ belongs to the Sobolev space of order $\alpha$, $H^\al(\Tor,\Rnu^n)$ 
and that $F$ belongs to $H^{\al}(\Tor,\Rnu^n)$,
where $\al$ is  bigger than $\frac{n}2+2$. The symbol
$\Tor$ denotes the $n$-dimensional torus, namely
$\Tor = \underbrace{S^1 \x \dots \x S^1}_n$, and $S^1$ is a unit circle.
We extend the functions $F$ and $h$ to $\Rnu^n$ periodically, and
use the same symbols for the extended functions.
Consider the following system of forward-backward SDEs in $H^\al(\Tor,\Rnu^n)$:
\aaa{
\label{fbsde-main}
\begin{cases}
  dZ_s^{t,e} = Y_s^{t,e} ds + \sqrt{2\nu}\, dW_s , \\
  dY_s^{t,e} = - F(s,\fdot) \circ Z_s^{t,e}\,ds + \sqrt{2\nu}\, X_s^{t,e}  dW_s, \\
  Z^{t,e}_t = e; \; Y_T^{t,e} = h \circ Z_T^{t,e}
\end{cases}
}
where $e: \Tor\to\Tor$ is the identical map,  $W_s$ is an $n$-dimensional
Brownian motion and $\nu >0$.
 Let $\mc F_s = \sg(W_r, r \in [0,s])$.
 We would like to find an $\mc F_s$-adapted
 triple of
 stochastic processes $(Z_s^{t,e}, Y_s^{t,e}, X_s^{t,e})$
with values in 
$H^\al(\Tor,\Rnu^n) \x H^\al(\Tor,\Rnu^n) \x H^\al(\Tor,\Rnu^{n\x n})$
which solves \eqref{fbsde-main}.
 Note that the process $X_s^{t,e}$
 takes values in the space of linear operators $\mc L(\Rnu^n, H^\s(\Tor,\Rnu^n))$
 ($\approx H^\s(\Tor,\Rnu^{n\x n})$), i.e.
 \aaa{
 \label{process-x}
  X_s^{t,e} = \sum_{i=1}^n X^{i}_s\ox e_i
 }
 where the processes $X^{i}_s$ take values in $H^\s(\Tor,\Rnu^n)$,
 and $\{e_i\}_{i=1}^n$ is an orthonormal basis of $\Rnu^n$. 
 Define
 \aaa{
 \label{K}
 K = \sup_{\Tor} |\nab h| +  \sup_{[0,T]\x \Tor} |\nab F|
 }
where $T$ is the arbitrary fixed number that we used to obtain
the backward Burgers equation.
Our main result is the following.
\begin{thm}
\label{main}
Let $h\in H^\al(\Tor,\Rnu^n)$, $F(s,\fdot)\in H^{\al}(\Tor,\Rnu^n)$, $s\in [0,T]$,  be 
such that $\nab^l h$ and $\nab^{l} F(s,\fdot)$ 
are bounded for all $l\lt \al$ and for all $s\in [0,T]$.
Then, there exists a $T_0>0$  that depends only on $\al$
and $K$ defined by \eqref{K} and such that for every $T<T_0$ and
for every $t\in [0,T]$, there exists a unique $\mc F_s$-adapted solution 
$(Z^{t,e}_s,Y^{t,e}_s, X^{t,e}_s)$  
to \eqref{fbsde-main} on $[t,T]$ with values in 
$H^\al(\Tor,\Rnu^n) \x H^\al(\Tor,\Rnu^n) \x H^\al(\Tor,\Rnu^{n\x n})$.
Moreover, the function $y:[0,T]\x\Tor\to \Rnu^n$, $(t,\te) \mto Y^{t,e}_t(\te)$ 
is deterministic, and solves the problem \eqref{brgs}.
\end{thm}
First we prove the existence and uniqueness of an $\mc F_s$-adapted 
$H^\al(\Tor,\Rnu^n) \x H^\al(\Tor,\Rnu^n) \x H^\al(\Tor,\Rnu^{n\x n})$-valued
solution $(Z^{t,\xi}_s,Y^{t,\xi}_s, X^{t,\xi}_s)$ 
to the following problem:
\aaa{
\label{fwd-xi}
 &Z^{t,\xi}_s = \xi + \int_t^s Y^{t,\xi}_r dr + \sqrt{2\nu}\, (W_s - W_t),\\
 \label{bwd-xi}
 &Y_s^{t,\xi} = h \circ Z_T^{t,\xi} + \int_s^T F(r,\fdot)\circ Z_r^{t,\xi}\,dr
 - \sqrt{2\nu}\, \int_s^T X_r^{t,\xi} dW_r,
 }
where $\xi$ is an $\mc F_t$-measurable $H^\al(\Tor,\Rnu^n)$-valued random variable.
Without loss of generality we will assume that the derivatives
$\nab^l h$ and $\nab^{l} F(t,\fdot)$, $t\in [0,T]$, are everywhere defined.
In the following, we will identify $H^0(\Tor,\Rnu^n)$ and $L_2(\Tor,\Rnu^n)$. 
The proof of Theorem \ref{main} will be devided in several lemmas.
\begin{lem}
\label{lem1}
Let $\E\,\|\xi\|^p_{L_p(\Tor,\Rnu^n)}$ and 
$\E\,\|\nab^i \xi\|^p_{L_p(\Tor,\Rnu^{n^i})}$ be bounded
for all integers $i\lt l$ and $p\gt 2$. Further suppose that the
FBSDEs \eqref{fwd-xi}, \eqref{bwd-xi} have a solution in $H^l(\Tor,\Rnu^n)$. 
Then for any integer $p\gt 2$ there exists a $T_0>0$ that depends only on $p$ 
and $K$ defined by \eqref{K} and
such that for all positive $T<T_0$, 
for all $i\lt l$, $\E\, \|Z_s^{t,\xi}\|^p_{L_p(\Tor,\Rnu^n)}$
and $\E\, \|\nab^i Z_s^{t,\xi}\|^p_{L_p(\Tor,\Rnu^{n^i})}$ are bounded on $[t,T]$.
\end{lem}
\begin{proof}
Everywhere below, for convenience, we use the same symbols $\gm$ and $\gm_i$, $i=1,2,3 \ldots$, for
(different) constants in different formulas. All $\gm$ and $\gm_i$ 
below are positive and do not depend on $s\in [t,T]$ and $\te\in\Tor$.
Note that, for any Hilbert norm,
\aa{
&\bigl(\|g\|^p\bigr)'h  = p\|g\|^{p-2}(g,h),\\
&\bigl(\|g\|^p\bigr)''h_1 h_2 = p(p-2)\|g\|^{p-4}(g,h_1)(g,h_2)+
p\|g\|^{p-2}(h_1,h_2).
}
Fix a $\te \in \Tor$, and let $z_s = Z^{t,\xi}_s(\te)$, $y_s = Y^{t,\xi}_s(\te)$, $x_s = X^{t,\xi}_s(\te)$.
BSDE \eqref{bwd-xi} and It\^o's formula imply:
\mm{
\E\,|y_s|^p 
- \sqrt{2\nu}\, p(p-2)\int_s^T \E\, \bigl[|y_r|^{p-4}\sum_{i=1}^n|(x^{i}_r,y_r)|^2\bigr]\,dr \\
-\sqrt{2\nu}\, p\int_s^T \hspace{-1mm} \E\,\bigl[|y_r|^{p-2}|x_r|^2\bigr]\,dr  = \E\,|h(z_T)|^p 
+ 2p\int_s^T \hspace{-1mm}\E\,\bigl[|y_r|^{p-2}(F(r,z_r), y_r)\bigr]\, dr
}
where $x^i_r = X^i_r(\te)$ and the processes $X^i_r$ were introduced in representation \eqref{process-x}. 
Taking into account that $h$ and $F$ are bounded on $\Tor$ and resp. $[0,T]\x\Tor$, and applying
Young's inequality we obtain the existence of constants 
$\gm_1$ and $\gm_2$ such that
\aa{
\E\, |y_s|^p \lt \gm_1 + 
\gm_2 \int_s^T \E\,|y_r|^{p}\, dr.
}
Applying Gronwall's lemma and then integrating over $\Tor$ we obtain
that there exists a constant $\gm$ such that
\aa{
\E\, \|Y^{t,\xi}_s\|^p_{L_p(\Tor,\Rnu^n)} \lt \gm.
}
From SDE \eqref{fwd-xi}, we deduce the existence of 
constants $\gm_1$ and $\gm_2$  such that
\aa{
\E\, |z_s|^p \lt \gm_1 \,\E\,|\xi|^p + \gm_2 \int_t^s \E\, |y_r|^p \, dr. 
}
Integrating  over $\Tor$ and modyfing $\gm$ we obtain:
\aa{
\E\, \|Z^{t,\xi}_s\|_{L_p(\Tor,\Rnu^n)}^p \lt \gm.
}
Let us prove now that 
$\E\, \|\nab Z^{t,\xi}_s\|^p_{L_p(\Tor,\Rnu^n)}$ is bounded, where $\nab$ stands, as usual, for the space derivative.
The triple $(\nab Z^{t,\xi}_s, \nab Y^{t,\xi}_s, \nab X^{t,\xi}_s)$ solves the FBSDEs:
 \aaa{
 \label{fwd-d1}
&  \nab Z_s^{t,\xi} = \nab \xi + \int_t^s \nab Y_r^{t,\xi} \, dr \\
 \label{bwd-d1}
& \nab Y_s^{t,\xi} = \nab h\bigl((Z_T^{t,\xi}(\fdot)\bigr)\nab Z_T^{t,\xi} + 
\int_s^T \hspace{-2mm}\nab F\bigl(r,Z_r^{t,\xi}(\fdot)\bigr)\nab Z_r^{t,\xi} \, dr 
- \sqrt{2\nu}\int_s^T\hspace{-2mm} \nab X_r^{t,\xi} \,dW_r.
 }
 For simplicity of the notation, let us introduce the processes  $ z_s = Z^{t,\xi}_s(\te)$, $y_s =  Y^{t,\xi}_s(\te)$ and $ x_s =  X^{t,\xi}_s(\te)$.
  It\^o's formula together with the BSDE \eqref{bwd-d1} imply:
  \mm{
|\nab y_s|^p 
+ \sqrt{2\nu}\, p(p-2)\int_s^T |\nab y_r|^{p-4}\sum_{i=1}^n|(\nab x^{i}_r,\nab y_r)|^2\,dr \\
+\sqrt{2\nu}\, p\int_s^T  |\nab y_s|^{p-2}|\nab x_r|^2\,dr  = |\nab h(z_T)\nab z_T|^p \\
+ 2p\int_s^T |\nab y_r|^{p-2}(\nab F(r,z_r)\nab z_r,\nab y_r)\,dr 
+ 2p\int_s^T |\nab y_r|^{p-2}(\nab y_r, \nab x_r\, dW_r).
}
This implies  that there exist constants $\gm_1$, $\gm_2$, and $\gm_3$ such that for all $s\in [t,T]$
\aa{
\E|\nab y_s|^p \lt \gm_1\,\E\, |\nab \xi|^p  
+ \gm_2 \,  \E
\int_t^T |\nab y_r|^{p-2}\,|\nab z_r|^2 \, dr
+\gm_3\,\int_t^T \hspace{-1mm} \E\, |\nab y_r|^p\, dr.
}
From Young's inequality,
\aaa{
\label{young_holder}
|\nab y_r|^{p-2}\,|\nab z_r|^2 \lt \frac{p-2}{p}|\nab y_r|^{p} + \frac{2}p\, |\nab z_r|^p.
}
Therefore, we can find constants $\gm_1$ and $\gm_2$ 
such that
\aa{
\E |\nab y_s|^p \lt \gm_1 \, \E\,|\nab \xi|^p + \gm_2\int_t^T \E |\nab y_r|^p \, dr.
}
Choosing $T_0$ smaller than $\frac1{\gm_2}$ we deduce that there exists a constant $\gm$
such that
\aa{
\E |\nab y_s|^p \lt \gm\,\E\, |\nab \xi|^p.
}
Integrating over $\Tor$ gives:
\aaa{
\label{est-dy}
\E \|\nab Y^{t,\xi}_s\|^p_{L_p(\Tor,\Rnu^{n^2})} \lt \gm\,\E\,\|\nab \xi\|^p_{L_p(\Tor,\Rnu^{n^2})}.
}
Next,  the SDE \eqref{fwd-d1} implies that there exist positive constants $\gm_3$ and $\gm_4$ such that
\aa{
\E\, \|\nab Z^{t,\xi}_s\|_{L_p(\Tor,\Rnu^{n^2})}^p \lt \gm_3\,\E\, \|\nab \xi\|_{L_p(\Tor,\Rnu^{n^2})}^p
+ \gm_4 \int_t^s \E\,\|\nab Y^{t,\xi}_r\|_{L_p(\Tor,\Rnu^{n^2})}^p \, dr.
}
Combining this and \eqref{est-dy} and modifying $\gm$ and $\gm_1$ we obtain:
\aa{
\E\, \|\nab Z^{t,\xi}_s\|_{L_p(\Tor,\Rnu^{n^2})}^p \lt \gm_1\,\E\, \|\nab \xi\|_{L_p(\Tor,\Rnu^{n^2})}^p \lt \gm.
}
Let us assume $\E\, \|\nab^i Z^{t,\xi}_s\|_{L_p(\Tor,\Rnu^{n^2})}^p \lt \gm$
for all integers $i\lt l-1$, and prove that 
$\E\, \|\nab^l Z^{t,\xi}_s\|_{L_p(\Tor,\Rnu^{n^2})}^p \lt \gm$The triple $(\nab^l Z^{t,\xi}_s, \nab^l Y^{t,\xi}_s, \nab^l X^{t,\xi}_s)$ solves the FBSDEs \eqref{fwd-dl},
\eqref{bwd-dl} below which are obtain from \eqref{fwd-xi}, \eqref{bwd-xi} by differentiating
both parts $l$ times:
 \aaa{
  \label{fwd-dl}
&\nab^l Z_s^{t,\xi} = \nab^l \xi + \int_t^s \nab^l Y_r^{t,\xi} \, dr \\
&\nab^l Y_s^{t,\xi} = \nab h\bigl(Z_T^{t,\xi}(\fdot)\bigr)\nab^l Z_T^{t,\xi} + \int_s^T \nab F\bigl(r,Z_r^{t,\xi}(\fdot)\bigr)\nab^l Z_r^{t,\xi} \, dr 
 \notag \\
 & +\sum_{j=2}^l \nab^j h\bigl(Z_T^{t,\xi}(\fdot)\bigr)\Bigl[\sum_{i_1+\dots+i_j = l-j+1}\nab^{i_1}Z_T^{t,\xi}\ldots \nab^{i_j}Z_T^{t,\xi}\Bigr]
 \notag \\
 & +\int_s^T \sum_{j=2}^{l} \nab^{j} F\bigl(r, Z_r^{t,\xi}(\fdot)\bigr)\Bigl[\sum_{i_1+\dots+i_j = l-j+1}\nab^{i_1}Z_r^{t,\xi}\ldots \nab^{i_j}Z_r^{t,\xi}\Bigr]\, dr
 \notag \\
  \label{bwd-dl}
&  -\sqrt{2\nu}\int_s^T \nab^l X_r^{t,\xi} \,dW_r.
 }
The argument below is similar to the one we have used for the first order derivatives.
It\^o's formula and the BSDE \eqref{bwd-dl} imply:
\mmm{
\label{ch4}
|\nab^l y_s|^p +\sqrt{2\nu}\, p\int_s^T  |\nab^l y_s|^{p-2}|\nab^l x_r|^2\,dr  \\
+ \sqrt{2\nu}\, p(p-2)\int_s^T |\nab^l y_r|^{p-4}\sum_{i=1}^n|(\nab^l x^{i}_r,\nab^l y_r)|^2\,dr \\
= \Bigl|\nab h(z_T)\nab^l z_T
+\sum_{j=2}^l \nab^j h(z_T)\Bigl[\sum_{i_1+\dots+i_j = l-j+1}\nab^{i_1}z_T\ldots \nab^{i_j}z_T\Bigr]
\Bigr|^p\\ 
+ 2p\int_s^T |\nab^l y_r|^{p-2}(\nab F(r,z_r)\nab^l z_r,\nab^l y_r)\,dr 
+ 2p\int_s^T |\nab^l y_r|^{p-2}(\nab^l y_r, \nab^l x_r\, dW_r)\\
  +\int_s^T |\nab^l y_r|^{p-2}\Bigl(\sum_{j=2}^{l} \nab^{j} F(r, z_r)\Bigl[\sum_{i_1+\dots+i_j = l-j+1}\nab^{i_1}z_r
  \ldots \nab^{i_j}z_r\Bigr],\nab^l y_r\Bigr)\, dr.
}
Note that by \eqref{fwd-dl}, there exist constants $\gm_1$ and $\gm_2$ so that
\aa{
|\nab^l z_s|^p \lt \gm_1 |\nab^l \xi|^p + \gm_2 \int_t^s |\nab^l y_r|^p\, dr
}
This and \eqref{ch4} imply for all $s\in [t,T]$ there exist 
constants $\gm_3$, $\gm_4$, $\gm_5$, $\gm_6$, and $\gm_7$
such that
\mm{
\E|\nab^l y_s|^p \lt \gm_3\,\E\, |\nab^l \xi|^p  
+ \gm_4 \,  \E 
\int_t^T |\nab^l y_r|^{p-2}\,|\nab^l z_r|^2 \, dr
+\gm_5\,\int_t^T \hspace{-1mm} \E\, |\nab^l y_r|^p\, dr\\
+ \gm_6\hspace{-2mm} 
\sum_{i_1+\dots+i_j \lt l-1}\E\, |\nab^{i_1}z_T\ldots \nab^{i_j}z_T|^p\\
+ \gm_7 \,\E\int_t^T \hspace{-1mm}  |\nab^l y_r|^{p-2}\, 
\Bigl|\sum_{j=2}^l \nab^{j} F(r, z_r)\,\Bigl[\sum_{i_1+\dots+i_j = l-j+1}\nab^{i_1}z_r
  \ldots \nab^{i_j}z_r\Bigr]\Bigr|^2\, dr
}
Note that we can 
apply inequality \eqref{young_holder} where $\nab y_r$ is replaced by $\nab^l y_r$ and $\nab z_r$ 
is replaced by $\nab^l z_r$, $r\in [t,T]$. 
Also, Young's inequality implies that 
there exists a constant $\gm_8$ such that
\mm{
 |\nab^l y_r|^{p-2}\, 
\Bigl|\sum_{j=2}^l \nab^{j} F(r, z_r)\,\Bigl[\sum_{i_1+\dots+i_j = l-j+1}\nab^{i_1}z_r
  \ldots \nab^{i_j}z_r\Bigr]\Bigr|^2 \\
  \lt \frac{p-2}{p}\, |\nab^l y_r|^p + 
  \frac2{p}\,\Bigl|\sum_{j=2}^l \nab^{j} F(r, z_r)\,\Bigl[\sum_{i_1+\dots+i_j = l-j+1}\nab^{i_1}z_r
  \ldots \nab^{i_j}z_r\Bigr]\Bigr|^p \\
  \lt \frac{p-2}{p}\, |\nab^l y_r|^p +
\gm_8 \sum_{i_1+\dots+i_j \lt l-1}|\nab^{i_1}z_r \dots \nab^{i_j}z_r|^p.
}
Finally we obtain that
there exist constants $\gm_1$, $\gm_2$, $\gm_3$, and $\gm_4$ such that 
\mm{
\E\,|\nab^l y_s|^p \lt \gm_1\, \E\,|\nab^l \xi|^p  +  \gm_2 \sum_{i_1+\dots+i_j \lt l-1}\E\, |\nab^{i_1}z_T\ldots \nab^{i_j}z_T|^p\\
+\gm_3 \int_t^T \sum_{i_1+\dots+i_j \lt l-1}\E\,|\nab^{i_1}z_r \dots \nab^{i_j}z_r|^p \, dr + \gm_4 \int_t^T \E\,|\nab^l y_r|^p \, dr.
}
Choosing $T_0$ smaller than $\frac1{\gm_4}$ and modifying $\gm_1$, $\gm_2$ and $\gm_3$ we obtain that
\mmm{
\label{es-yl}
\E\,|\nab^l y_s|^p \lt \gm_1\, \E\, |\nab^l \xi|^p  +  \gm_2 \sum_{i_1+\dots+i_j \lt l-1}\E\, |\nab^{i_1}z_T\ldots \nab^{i_j}z_T|^p\\
+\gm_3 \int_t^T \sum_{i_1+\dots+i_j \lt l-1}\E\,|\nab^{i_1}z_r \dots \nab^{i_j}z_r|^p \, dr
}
and moreover, modifying $\gm_1$, $\gm_2$ and $\gm_3$ we obtain that
\mmm{
\label{est-zl}
\E|\nab^l z_s|^p \lt \gm_1\, \E\, |\nab^l \xi|^p  +  \gm_2 \sum_{i_1+\dots+i_j \lt l-1}\E\, |\nab^{i_1}z_T\ldots \nab^{i_j}z_T|^p\\
+\gm_3 \int_t^T \sum_{i_1+\dots+i_j \lt l-1}\E\,|\nab^{i_1}z_r \dots \nab^{i_j}z_r|^p \, dr.
}
Integrating over $\Tor$ and taking into account that $\E\,\|\nab^l \xi\|^p_{L_p(\Tor,\Rnu^{n^l})}$ is bounded by assumption,
 $\E\,\|\nab^{i_1}Z^{t,\xi}_r\ldots \nab^{i_j}Z^{t,\xi}_r\|^p_{L_p(\Tor,\Rnu^{n^l})}$, $r\in [t,T]$, are bounded by 
H\"older's inequality and the induction hypothesis, 
we obtain that $\E\,\|\nab^l Y^{t,\xi}_s\|^p_{L_p(\Tor,\Rnu^{n^l})}$ and $\E\,\|\nab^l Z^{t,\xi}_s\|^p_{L_p(\Tor,\Rnu^{n^l})}$
are bounded.
\end{proof}
\begin{lem}
\label{lemt0}
There exists a $T_0 > 0$ such that for every positive $T<T_0$ and for every $t\in [0,T]$, FBSDEs \eqref{fwd-xi}, \eqref{bwd-xi}
has a unique $\mc F_s$-adapted solution on $[t,T]$ with values in
$H^\al(\Tor,\Rnu^n)\x H^\al(\Tor,\Rnu^n)\x H^\al(\Tor,\Rnu^{n\x n})$.
\end{lem}
\begin{proof}
First we prove the existence of solution in $L_2(\Tor,\Rnu^n)$.
Let us consider the map
\aa{
\Gamma: L_2(\Om,L_2(\Tor,\Rnu^n)) \to L_2(\Om,L_2(\Tor,\Rnu^n)), \; Y^{t,\xi}_s \to \bar Y^{t,\xi}_s
}
which is defined by the FBSDEs below:
\aaa{
\label{1sde}
&\bar Z^{t,\xi}_s = \xi+ \int_t^s Y^{t,\xi}_r \, dr + \sqrt{2\nu}\, (W_s - W_t), \\
\label{2sde}
&\bar Y^{t,\xi}_s = h\bigl(\bar Z^{t,\xi}_T(\fdot)\bigr)  + \int_s^T F \bigl(r,\bar Z^{t,\xi}_r(\fdot)\bigr) \, dr
- \sqrt{2\nu}\int_s^T \bar X^{t,\xi}_r \, dW_r.
}
First we find $\bar Z^{t,\xi}_s$ from the SDE \eqref{1sde}, and substitute it into BSDE \eqref{2sde}. 
Then we find the unique $\mc F_s$-adapted solution $(\bar Y^{t,\xi}_s,\bar X^{t,\xi}_s)$
of BSDE \eqref{2sde}.
Namely,
\aa{
\bar Y^{t,\xi}_s = \E \, [h\bigl(\bar Z^{t,\xi}_T(\fdot)\bigr) + \int_s^T F\bigl(r,\bar Z^{t,\xi}_r(\fdot)\bigr)\, dr \,|\, \mc F_s ],
}
and $\bar X^{t,\xi}_s$ exists by the martingale representation theorem (see \cite{ns1}).
Note that since $h$ and $F$ are bounded on $\Tor$ and resp. on  $[0,T]\x \Tor$, $\bar Y^{t,\xi}_s$ takes values
in $L_2(\Tor,\Rnu^n)$.
The process $\bar X^{t,\xi}_s$ is actually not needed for the definition of the map $\Gamma$.
Let us prove that the map $\Gamma$ is a contraction. 
Let $V^{t,\xi}_s$, $Y^{t,\xi}_s$ $\in L_2(\Om, L_2(\Tor,\Rnu^n))$, and 
$\bar V^{t,\xi}_s = \Gamma(V^{t,\xi}_s)$, $\bar Y^{t,\xi}_s = \Gamma(Y^{t,\xi}_s)$.
Further let $\bar U^{t,\xi}_s$ be obtained from \eqref{1sde}, i.e.
$\bar U^{t,\xi}_s = \xi+ \int_t^s V^{t,\xi}_r \, dr + \sqrt{2\nu}\, (W_s - W_t)$,
and $\bar W^{t,\xi}_s$ be obtained from the  BSDE \eqref{2sde}, i.e.  
$
\bar V^{t,\xi}_s = h\bigl(\bar U^{t,\xi}_T(\fdot)\bigr)  + \int_s^T F\bigl(r,\bar U^{t,\xi}_r(\fdot)\bigr) \, dr
- \sqrt{2\nu} \int_s^T \bar W^{t,\xi}_r \, dW_r.
$
The SDE \eqref{1sde} implies that for all $s\in [t,T]$ 
\aaa{
\label{1ineqv}
\|\bar Z^{t,\xi}_s - \bar U^{t,\xi}_s\|^2_{L_2(\Tor,\Rnu^n)} \lt 
(s-t) \int_t^s \|Y^{t,\xi}_r - V^{t,\xi}_r\|^2_{L_2(\Tor,\Rnu^n)} \, dr.
}
The SDE \eqref{2sde} and It\^o's formula imply that 
\aa{
&\E\,\|\bar Y^{t,\xi}_s - \bar V^{t,\xi}_s\|^2_{L_2(\Tor,\Rnu^n)}
+2\nu \int_s^T \E\, \|\bar X^{t,\xi}_r - \bar W^{t,\xi}_r\|^2_{L_2(\Tor,\Rnu^n)} dr\\
&= \E\,\|h\bigl(\bar Z^{t,\xi}_T(\fdot)\bigr)- h\bigl(\bar U^{t,\xi}_T(\fdot)\bigr)\|^2_{L_2(\Tor,\Rnu^n)} \\
&+ 2 \int_s^T \E\, ( F\bigl(r, \bar Z^{t,\xi}_r(\fdot)\bigr) 
- F\bigl(r,\bar U^{t,\xi}_T(\fdot)\bigr),\bar Y^{t,\xi}_r - \bar V^{t,\xi}_r)_{L_2(\Tor,\Rnu^n)}\, dr.
}
Hence, 
\mm{
\E\,\|\bar Y^{t,\xi}_s - \bar V^{t,\xi}_s\|^2_{L_2(\Tor,\Rnu^n)} \lt 
\E\,\|h\bigl(\bar Z^{t,\xi}_T(\fdot)\bigr)- h\bigl(\bar U^{t,\xi}_T(\fdot)\bigr)\|^2_{L_2(\Tor,\Rnu^n)} \\
+ \int_t^T \E \,\|F\bigl(r,\bar Z^{t,\xi}_r(\fdot)\bigr) - F\bigl(r,\bar U^{t,\xi}_r(\fdot)\bigr)\|^2_{L_2(\Tor,\Rnu^n)} dr
\\+ \int_s^T \E\,\|\bar Y^{t,\xi}_r - \bar V^{t,\xi}_r\|^2_{L_2(\Tor,\Rnu^n)}dr. 
}
Gronwall's lemma and inequality \eqref{1ineqv} 
imply that 
\aa{
\E\,\|\bar Y^{t,\xi}_s - \bar V^{t,\xi}_s\|^2_{L_2(\Tor,\Rnu^n)} \lt 
\td\gm(T) \int_t^T \E\,\|Y^{t,\xi}_r - V^{t,\xi}_r\|^2_{L_2(\Tor,\Rnu^n)} dr
}
where 
$\td\gm(T) = e^T K$ and $K$ is defined by \eqref{K}. 
Let us pick  $T_0$ so that $\gm(T_0) =T_0\td\gm(T_0) < 1$.
Then, if $T<T_0$, 
\aaa{
\label{1fp-est}
\sup_{s\in [t,T]}\E\,\|\bar Y^{t,\xi}_s - \bar V^{t,\xi}_s\|_{L_2(\Tor,\Rnu^n)}
\lt \gm(T) \sup_{s\in [t,T]}\E\,\|Y^{t,\xi}_s - V^{t,\xi}_s\|_{L_2(\Tor,\Rnu^n)}.
}
This proves that for $T<T_0$ there is a unique fixed point $Y^{t,\xi}_s$ of the map $\Gamma$.
The processes $Z^{t,\xi}_s$ and  $X^{t,\xi}_s$ can be determined  from \eqref{bwd-xi}
as described above.
Let us consider now the FBSDEs which is obtained from \eqref{fwd-xi}, \eqref{bwd-xi} by differentiation
with respect to $\te\in\Tor$:
 \aaa{
 \label{fbsde-d1}
 \begin{cases}
  \nab Z_s^{t,\xi} = \nab \xi + \int_t^s \nab Y_r^{t,\xi} \, dr \\
 \nab Y_s^{t,\xi} = \nab h\bigl(Z_T^{t,\xi}(\fdot)\bigr)\nab Z_T^{t,\xi} + 
  \int_s^T \nab F\bigl(r,Z_r^{t,\xi}(\fdot)\bigr)\nab Z_r^{t,\xi} \, dr \\
- \sqrt{2\nu}\int_s^T \nab X_r^{t,\xi} \,dW_r.
 \end{cases}
 }
 Now we assume that the solution $(Z^{t,\xi}_s, Y^{t,\xi}_s, X^{t,\xi}_s)$ is known, and therefore the FBSDEs \eqref{fbsde-d1}
 are regarded as a system of SDEs with random coefficients. 
 Clearly, if we prove the existence of solution $(\nab Z^{t,\xi}_s, \nab Y^{t,\xi}_s, \nab X^{t,\xi}_s)$
 to \eqref{fbsde-d1} it would imply that the solution $(Z^{t,\xi}_s, Y^{t,\xi}_s, X^{t,\xi}_s)$
 is differentiable in $\te$, and solves \eqref{fwd-xi}, \eqref{bwd-xi} in $H^1(\Tor,\Rnu^n)$.
 The proof of this fact uses standard approaches described for example in \cite{delarue} or
 \cite{belopolskaya}.
 The same argument as before applied to the triple $(\nab Z^{t,\xi}_s, \nab Y^{t,\xi}_s, \nab X^{t,\xi}_s)$
 as well as the boundedness of $\nab h$ and $\nab F$ on $\Tor$
 and resp. on $[0,T]\x \Tor$
 imply the existence and uniqueness of a solution to \eqref{fbsde-d1}, and therefore
 the existence and uniqueness of a solution to \eqref{fwd-xi}, \eqref{bwd-xi} 
 with respect to the $H^1(\Tor,\Rnu^n)$-norm.
 Indeed, consider the map 
 \aa{
 \Gamma^{(1)}: L_2(\Om,L_2(\Tor,\Rnu^{n^2})) \to L_2(\Om,L_2(\Tor,\Rnu^{n^2})), 
 \; \nab Y^{t,\xi}_s \to \ovl{\nab Y^{t,\xi}_s}
 }
 which is defined by the FBSDEs:
 \aaa{
 \label{sde1-d1}
 &\ovl{\nab Z^{t,\xi}_s} = \nab \xi + \int_t^s \nab Y^{t,\xi}_r dr,\\
 \label{sde2-d1}
 &\ovl{\nab Y^{t,\xi}_s} = \nab h\bigl(Z_T^{t,\xi}(\fdot)\bigr)\ovl{\nab Z^{t,\xi}_T} + \int_s^T \nab 
 F\bigl(r,Z^{t,\xi}_r(\fdot)\bigr)\ovl{\nab Z^{t,\xi}_r} \, dr \\
&-\sqrt{2\nu} \int_s^T \ovl{\nab X^{t,\xi}_r}\, dW_r. \notag
 }
The process $\ovl{\nab Z^{t,\xi}_s}$ is obtained from  \eqref{sde1-d1}, and the processes
$\ovl{\nab Y^{t,\xi}_s}$ and $\ovl{\nab X^{t,\xi}_s}$ are obtained from the BSDE \eqref{sde2-d1}
as its unique $\mc F_s$-adapted solution.
Let $\nab V^{t,\xi}_s$, $\ovl{ \nab U^{t,\xi}_s}$, $\ovl{ \nab V^{t,\xi}_s}$, and $\ovl{ \nab W^{t,\xi}_s}$
be associated with the map $\Gamma^{(1)}$ and correspond to the processes
$V^{t,\xi}_s$, $\bar U^{t,\xi}_s$, $\bar V^{t,\xi}_s$, and $\bar W^{t,\xi}_s$
in the fixed point argument for the map $\Gamma$.
The SDE \eqref{sde1-d1} implies the estimate:
\aa{
\|\ovl{\nab Z^{t,\xi}_s}-\ovl{\nab U^{t,\xi}_s}\|^2_{L_2(\Tor,\Rnu^{n^2})} \lt (s-t) 
\int_t^s\|{\nab Y^{t,\xi}_r}-{\nab V^{t,\xi}_r}\|^2_{L_2(\Tor,\Rnu^{n^2})} dr.
 }
Application of It\^o's formula to
$\|\ovl{\nab Y^{t,\xi}_r}-\ovl{\nab V^{t,\xi}_r}\|^2_{L_2(\Tor,\Rnu^{n^2})}$ 
gives
\aa{
&\E\,\|\ovl{\nab Y^{t,\xi}_s}-\ovl{\nab V^{t,\xi}_s}\|^2_{L_2(\Tor,\Rnu^{n^2})}
+2\nu\int_s^T\E \, \|\ovl{\nab X^{t,\xi}_r}-\ovl{\nab W^{t,\xi}_r}\|^2_{L_2(\Tor,\Rnu^{n^3})}\\
&= \E\,\|\nab h\bigl(Z^{t,\xi}_T(\fdot)\bigr)(\ovl{\nab Z^{t,\xi}_T} - \ovl{\nab U^{t,\xi}_T})\|^2_{L_2(\Tor,\Rnu^{n^2})}\\
&+2 \int_s^T \E\,\bigl(\nab F\bigl(r,Z^{t,\xi}_r(\fdot)\bigr)(\ovl{\nab Z^{t,\xi}_r} - \ovl{\nab U^{t,\xi}_r}),
\ovl{\nab Y^{t,\xi}_r}-\ovl{\nab V^{t,\xi}_r}\bigr)_{L_2(\Tor,\Rnu^{n^2})}.
}
The same argument as for the map $\Gamma$ implies that
\aa{
\sup_{s\in [t,T]}\E\,\|\ovl{\nab Y^{t,\xi}_s}-\ovl{\nab V^{t,\xi}_s}\|^2_{L_2(\Tor,\Rnu^{n^2})} \lt \gm(T) 
\sup_{s\in [t,T]}\E\,\|{\nab Y^{t,\xi}_s}-{\nab V^{t,\xi}_s}\|^2_{L_2(\Tor,\Rnu^{n^2})} 
}
where $\gm(T)$ and $T$ can be choosen in exactly the same as in \eqref{1fp-est}.
Now let us assume that we proved the existence of solution to \eqref{fwd-xi}, \eqref{bwd-xi}
in $H^{l-1}(\Tor,\Rnu^n)$. Namely, we formally differentiate \eqref{fwd-xi}, \eqref{bwd-xi}
$l-1$ times with respect to $\te$,
and assume that we have proved the existence of a solution
$(\nab^{l-1} Z^{t,\xi}_s, \nab^{l-1} Y^{t,\xi}_s, \nab^{l-1} X^{t,\xi}_s)$
in the space $L_2(\Tor,\Rnu^{n^{l-1}})$.
Let us differentiate the FBSDE \eqref{fwd-xi}, \eqref{bwd-xi} once again. We obtain the FBSDEs 
\eqref{fwd-dl}, \eqref{bwd-dl} which we consider as the FBSDEs 
in $L_2(\Tor,\Rnu^{n^{l}})$ with random coefficients with respect to three unknown processes 
$(\nab^{l} Z^{t,\xi}_s, \nab^{l} Y^{t,\xi}_s, \nab^{l} X^{t,\xi}_s)$.
Consider the map
\aa{
 \Gamma^{(l)}: L_2(\Om,L_2(\Tor,\Rnu^{n^l})) \to L_2(\Om,L_2(\Tor,\Rnu^{n^l})), 
 \; \nab^l Y^{t,\xi}_s \to \ovl{\nab^l Y^{t,\xi}_s}
 }
 which is defined as follows: first we determine $\ovl{\nab^l Z^{t,\xi}_s}$ from the SDE
 \aa{
 \ovl{\nab^l Z^{t,\xi}_s} = \nab^l \xi + \int_t^s \nab^l Y^{t,\xi}_r dr.
 }
 Then we substitute $\ovl{\nab^l Z^{t,\xi}_s}$ in the SDE
 \aaa{
 &\ovl{\nab^l Y^{t,\xi}_s} = \nab h\bigl(Z_T^{t,\xi}(\fdot)\bigr)\ovl{\nab^l Z^{t,\xi}_T} + \int_s^T \nab 
 F\bigl(r,Z^{t,\xi}_r(\fdot)\bigr)\ovl{\nab^l Z^{t,\xi}_r} \, dr\notag\\  &+\sqrt{2\nu} \int_s^T \ovl{\nab^l X^{t,\xi}_r}\, dW_r
   +\sum_{j=2}^l \nab^j h\bigl(Z_T^{t,\xi}(\fdot)\bigr)\Bigl[\sum_{i_1+\dots+i_j = l-j+1}\nab^{i_1}Z_T^{t,\xi}\ldots \nab^{i_j}Z_T^{t,\xi}\Bigr]
 \notag \\
\label{sde2-dl}
 & +\int_s^T \sum_{j=2}^{l} \nab^{j} F\bigl(r, Z_r^{t,\xi}(\fdot)\bigr)
 \Bigl[\sum_{i_1+\dots+i_j = l-j+1}\nab^{i_1}Z_r^{t,\xi}\ldots \nab^{i_j}Z_r^{t,\xi}\Bigr]\, dr
}
and find a couple $\bigl(\ovl{\nab^l Y^{t,\xi}_s}, \ovl{\nab^l X^{t,\xi}_s}\bigr)$ as the unique $\mc F_s$-adapted
solution to \eqref{sde2-dl}.
Namely, we have the following expression for $\ovl{\nab^l Y^{t,\xi}_s}$:
\aaa{
 &\ovl{\nab^l Y^{t,\xi}_s} = \E\,\Bigl[\nab h\bigl(Z_T^{t,\xi}(\fdot)\bigr)\ovl{\nab^l Z^{t,\xi}_T} + \int_s^T \nab 
 F\bigl(r,Z^{t,\xi}_r(\fdot)\bigr)\ovl{\nab^l Z^{t,\xi}_r} \, dr \notag\\
  & +\sum_{j=2}^l \nab^j h\bigl(Z_T^{t,\xi}(\fdot)\bigr)\Bigl[\sum_{i_1+\dots+i_j = l-j+1}\nab^{i_1}Z_T^{t,\xi}\ldots \nab^{i_j}Z_T^{t,\xi}\Bigr]
\notag \\
\label{sde1-dl-1}
 & +\int_s^T \sum_{j=2}^{l} \nab^{j} F\bigl(r, Z_r^{t,\xi}(\fdot)\bigr)\Bigl[\sum_{i_1+\dots+i_j = l-j+1}\nab^{i_1}Z_r^{t,\xi}\ldots \nab^{i_j}Z_r^{t,\xi}\Bigr]\, dr
\, |\,\mc F_s\Bigr].
 }
 By Lemma \ref{lem1}, the last two terms of \eqref{sde1-dl-1} belong to $L_2(\Tor,\Rnu^{n^l})$, and therefore
$\ovl{\nab^l Y^{t,\xi}_s}$ takes values in $L_2(\Tor,\Rnu^{n^l})$. As before, we find
$\ovl{\nab^l X^{t,\xi}_r}$ by the martingale representation theorem.
 Since the coefficients of $\nab^l Z^{t,\xi}_T$ and of $\nab^l Y^{t,\xi}_r$ under the integral sign 
are the same as in the case $l=1$,
the fixed point argument will be also the same as for this case. In particular, 
$T_0$ will be the same as for the FBSDEs \eqref{fbsde-d1} and \eqref{fwd-xi}, \eqref{bwd-xi}.
By induction, we conclude that 
\eqref{fwd-dl}, \eqref{bwd-dl} has $\mc F_s$-adapted solutions 
$(\nab^{l} Z^{t,\xi}_s, \nab^{l} Y^{t,\xi}_s, \nab^{l} X^{t,\xi}_s)$
for every $l=1,\ldots,\al$. This implies that there exists an $\mc F_s$-adapted 
$H^\al(\Tor,\Rnu^n)\x H^\al(\Tor,\Rnu^n)\x H^\al(\Tor,\Rnu^{n\x n})$-solution to \eqref{fwd-xi}, \eqref{bwd-xi}.

Uniqueness of solution can be shown as in the proof of Lemma 15 of \cite{ns1}. 
\end{proof}

We have now shown the existence of solution for the forward-backward system of 
stochastic equations (2). From here we proceed to obtain the deterministic 
function $y$ which actually determines the drift of the process $Z_s^{t,e}$. 
This procedure is  the same that we have followed in [C-S] to derive the 
solution of Navier-Stokes equations from the solution of the corresponding 
stochastic system. The  difference is that, since now we are dealing with 
Burgers equation the incompressibility condition ($\div y=0$) is not present and, accordingly, the process $Z_s^{t,e}$ here belongs 
to the group $G^{\alpha}$ of $H^\al$-diffeomorphisms $\Tor\to\Tor$
whereas in [C-S] the corresponding relevant space is the subgroup of the volume-preserving diffeomorphisms.
Still, for sake of completeness, we present here the rest of the proof.
Everywhere below we assume that $T<T_0$, where $T_0$ is defined in Lemma \ref{lemt0}.
\begin{lem}
\label{lem2}
If $\xi$ is an $\mc F_t$-measurable $H^\al(\Tor,\Rnu^n)$-valued random variable,
then the solution $(Z_s^{t,\xi}, Y_s^{t,\xi}, X_s^{t,\xi})$ to 
\eqref{fwd-xi}, \eqref{bwd-xi} takes the form:
\aaa{
 \label{form_of_solution}
 (Z_s^{t,\xi}, Y_s^{t,\xi}, X_s^{t,\xi}) =   
 (Z^{t,e}_s\circ \xi,Y^{t,e}_s\circ \xi, X^{t,e}_s\circ \xi).
 }
\end{lem}
\begin{proof}
It suffices to prove the statement of the lemma in the space $L_2(\Tor,\Rnu^n)$.
Indeed, by uniqueness of solution,  $(Z^{t,\xi}_s, Y^{t,\xi}_s, X^{t,\xi}_s)$ 
is the unique solution to \eqref{fwd-xi}, \eqref{bwd-xi} in $L_2(\Tor,\Rnu^n)$, 
and therefore if we prove
\eqref{form_of_solution} in  $L_2(\Tor,\Rnu^n)$, then the triple
$(Z^{t,e}_s\circ \xi,Y^{t,e}_s\circ \xi, X^{t,e}_s\circ \xi)$ is the unique
solution to \eqref{fwd-xi}, \eqref{bwd-xi} also in $H^\al(\Tor,\Rnu^n)$.

Let us prove the statement in $L_2(\Tor,\Rnu^n)$.
We apply the operator $R_\xi$ of the composition with $\xi$ to  both parts of the
SDEs \eqref{fwd-e} and \eqref{bwd-e}:
\aaa{
\label{fwd-e}
 &Z^{t,e}_s = e + \int_t^s Y^{t,e}_r dr + \sqrt{2\nu}\, (W_s - W_t),\\
 \label{bwd-e}
 &Y_s^{t,e} = h \circ Z_T^{t,e} + \int_s^T F(r,\fdot)\circ Z_r^{t,e}\,dr
 - \sqrt{2\nu}\, \int_s^T X_r^{t,e} dW_r.
 }
Let us observe that we can write $R_\xi$ under the integral signs of the Bochner integrals. 
Indeed, since these integrals converge in $H^\al(\Tor,\Rnu^n)$, they also converge
with respect to at least the $\C(\Tor,\Rnu^n)$-topology. 
Due to the periodicity of
the functions under the (Bochner) integrals signs,
the composition of the integrands
with $\xi$ will preserve the convergence with respect to the $\C(\Tor,\Rnu^n)$-topology.
For the integrals converging in the $\C(\Tor,\Rnu^n)$-topology,
we can easily see that we are allowed to write $R_\xi$ under the integral signs.
This implies that if we consider convergence of the Bochner integrals with respect to
the $L_2(\Tor,\Rnu^n)$-topology, we can also write $R_\xi$ under the integral
signs.

Let us prove now that we are allowed to write
$R_\xi$ under the integral signs 
of the stochastic integrals in \eqref{fwd-e}, \eqref{bwd-e}.
First we prove this for the case when $\xi = g$ is deterministic.
Let $s$ and $S$ be such that
$t\lt s < S \lt T$, and let $\Phi_r$ be an $\mc F_r$-adapted
stochastically integrable process, and let $I(\Phi_r)$ denote
$\int_s^S \Phi_r \, dW_r$.
Let $\Phi^{(m)}_r$ be a sequence of simple stochastic processes
such that $I(\Phi^{(m)}_r)$ converge to 
$I(\Phi_r)$ with respect to the $L_2(\Om) \x H^\al(\Tor,\Rnu^n)$-norm.
Note that if $F \in H^\al(\Tor,\Rnu^n)$ extended to $\Rnu^n$ periodically, then
there exist constants $\gm_1$ and $\gm_2$ not depending on $F$ and such that
\aaa{
\label{chain1}
\|F\circ g\|_{L_2(\Tor,\Rnu^n)}\lt \gm_1 \|F\circ g\|_{\C(\Tor,\Rnu^n)} \lt 
\gm_1 \|F\|_{\C(\Tor,\Rnu^n)}\lt\gm_2 \|F\|_{H^\al(\Tor,\Rnu^n)}.
}
Therefore, since $\E\int_s^S \|\Phi^{(m)}_r - \Phi_r\|^2_{H^\al(\Tor,\Rnu^n)}dr \to 0$,
then $\E\int_s^S \|\Phi^{(m)}_r\circ g - \Phi_r \circ g\|^2_{L_2(\Tor,\Rnu^n)}dr\to 0$.
By It\^o's isometry, $\E \|I(\Phi^{(m)}_r\circ g) - I(\Phi_r \circ g)\|^2_{L_2(\Tor,\Rnu^n)}\to 0$. 
Again using \eqref{chain1}, we conclude that 
$\E \|I(\Phi^{(m)}_r)\circ g - I(\Phi_r)\circ g\|^2_{L_2(\Tor,\Rnu^n)}\to 0$
because 
$\E \|I(\Phi^{(m)}_r) - I(\Phi_r)\|^2_{H^\al(\Tor,\Rnu^n)}\to 0$.
Clearly, for simple stochastic processes it holds that $I(\Phi^{(m)}_r)\circ g = I(\Phi^{(m)}_r \circ  g)$,
and therefore $I(\Phi_r)\circ g = I(\Phi_r \circ  g)$.

Now let us take an $\mc F_t$-measurable stepwise function
$\xi=\sum_{i=1}^\infty g_i\ind_{A_i}$,
where $g_i\in H^\al(\Tor,\Rnu^n)$ and  the sets $A_i$ are $\mc F_t$-measurable.
We obtain:
\mm{
\int_s^S \Phi_r \, dW_r \circ \sum_{i=1}^\infty g_i \ind_{A_i} = 
\sum_{i=1}^\infty \ind_{A_i} \int_s^S \Phi_r \circ g_i \, dW_r =
\sum_{i=1}^\infty \int_s^S \ind_{A_i} \Phi_r \circ g_i \, dW_r \\
= \int_s^S \Phi_r \circ \sum_{i=1}^\infty g_i \ind_{A_i} \, dW_r.
}
Next,  we find a sequence of $\mc F_t$-measurable stepwise functions 
converging to $\xi$ in the space of continuous functions $\C(\Tor,\Rnu^n)$. 
This is possible due to
the separability of $\C(\Tor,\Rnu^n)$. Indeed, let us consider a countable
number of disjoint Borel sets $O^n_i$ covering $\C(\Tor,\Rnu^n)$,
and such that their diameter in the norm of $\C(\Tor,\Rnu^n)$ 
is smaller than $\frac1{n}$. 
Let $A_i^n=\xi^{-1}(O_i^n)$ and $g_i^n \in O_i^n$.
Define $\xi_n = \sum_{i=1}^\infty g_i^n \ind_{A_i^n}$. 
Then  for all $\om\in\Om$, we have $\|\xi-\xi_n\|_{\C(\Tor,\Rnu^2)} < \frac1{n}$.
As before, $I(\Phi)$ and $I(\Phi\circ \xi)$ denote
$\int_s^S \Phi_r\, dW_r$ and  $\int_s^S \Phi_r\circ \xi \, dW_r$ resp.
We have to prove that a.s. $I(\Phi)\circ \xi = I(\Phi\circ \xi)$.
For this it suffices to prove that
\aaa{
\label{61}
\lim_{n\to\infty} \E \|I(\Phi)\circ \xi_n - I(\Phi)\circ \xi\|^2_{L_2(\Tor,\Rnu^2)} &=0, \\
\label{62}
\lim_{n\to\infty} \E \|I(\Phi\circ \xi_n) - I(\Phi\circ \xi)\|^2_{L_2(\Tor,\Rnu^2)} &=0.
}
By \eqref{chain1}, $\|I(\Phi)\circ \xi_n\|_{L_2(\Tor,\Rnu^n)} \lt \gm_2 \|I(\Phi)\|_{H^\al(\Tor,\Rnu^n)}$,
and $\|I(\Phi)\circ \xi\|_{L_2(\Tor,\Rnu^n)} \lt \gm_2 \|I(\Phi)\|_{H^\al(\Tor,\Rnu^n)}$.
By Lebesgue's theorem, in \eqref{61} we can 
pass to the limit under the expectation sign.
 Relation \eqref{61} holds then by
the continuity of $I(\Phi)$ in $\te\in \Tor$. 
To prove \eqref{62} we observe that by It\^o's isometry, the limit in \eqref{62} equals to 
$\lim_{n\to\infty}\E\int_s^S \|\Phi_r\circ \xi_n - \Phi_r \circ \xi\|^2_{L_2(\Tor,\Rnu^2)} dr$.
The same argument that we used to prove \eqref{61} implies that we can pass  to
the limit under the expectation and the integral signs.
Relation \eqref{62} follows from the continuity of $\Phi_r$ in $\te\in\Tor$.
Hence,
 $(Z^{t,e}_s\circ \xi,Y^{t,e}_s\circ \xi, X^{t,e}_s\circ \xi)$
 is a solution to \eqref{fwd-xi}, \eqref{bwd-xi} in $L_2(\Tor,\Rnu^n)$. 
\end{proof}
\begin{lem}
\label{lem6}
The processes $Z^{t,e}_s$ and $Y^{t,e}_s$
have continuous path modifications. Namely, for these modifications
it holds that the
trajectories $[t,T]\to L_2(\Tor,\Rnu^n)$, $s\mto Z^{t,e}_s$
 and $[t,T]\to L_2(\Tor,\Rnu^n)$, $s\mto Y^{t,e}_s$ are continuous
with probability $1$.
\end{lem}
\begin{proof}
Let $s' > s$.
Application of It\^o's formula together with the BSDE \eqref{bwd-e} imply:
\aa{
\|Y^{t,e}_s - Y^{t,e}_{s'}\|^2_{L_2(\Tor,\Rnu^n)}
\lt \int_s^{s'}\hspace{-2mm} \|F(r, Z^{t,e}_r)\|^2_{L_2(\Tor,\Rnu^n)}\, dr + \int_s^{s'} 
\hspace{-2mm}\|Y^{t,e}_r - Y^{t,e}_{s'}\|^2_{L_2(\Tor,\Rnu^n)}\, dr.
}
Gronwall's lemma implies that there exist constants $\td\gm>0$ and $\gm>0$ such that
\aa{
\|Y^{t,e}_s - Y^{t,e}_{s'}\|^2_{L_2(\Tor,\Rnu^n)}
\lt\td \gm \int_s^{s'}\|F(r,Z^{t,e}_r)\|^2_{L_2(\Tor,\Rnu^n)}\, dr 
\lt \gm (s'-s).
}
This implies that if $p>1$ then
\aa{
\|Y^{t,e}_s - Y^{t,e}_{s'}\|^{2p}_{L_2(\Tor,\Rnu^n)}
\lt \gm |s-s'|^p.
}
By Kolmogorov's continuity criteria, $Y^{t,e}_s$ has a continuous path modification
with respect to the ${L_2(\Tor,\Rnu^n)}$-topology. The SDE \eqref{fwd-e} implies 
that $Z^{t,e}_s$ has a continuous path modification in the ${L_2(\Tor,\Rnu^n)}$-topology as well.
\end{proof}

Lemmas \ref{lem4} and \ref{lem7} below 
characterize the deterministic nature of the process $Y^{t,e}_s$ 
and describe its continuity properties.
\begin{lem}
\label{lem4}
The map
\aa{
[0,T]\x\Tor \to \Rnu^n, \; (t,\te) \to Y^{t,e}_t(\te) 
}
is deterministic and
the function $[0,T]\to H^\al(\Tor,\Rnu^n), t\mto Y^{t,e}_t$ is continuous.
\end{lem}
\begin{proof}
The first statement is a consequence of Blumenthal's zero-one law and 
the fact that the random variable $Y^{t,e}_t$ is $\mc F_0$-measurable
(as in Lemma 13 of [C-S] or Corollary 1.5. of [D]).

The proof of the continuity of the map
$[0,T]\to L_2(\Tor,\Rnu^n), t\mto Y^{t,e}_t$
follows as in Lemma 14 of \cite{ns1}.
Consider the FBSDEs below  on the interval $[0,T]$ with respect to 
$(\nab^l Z^{t,\xi}_s, \nab^l Y^{t,\xi}_s, \nab^l X^{t,\xi}_s)$:
 \aaa{
 \label{fwd-dl-ext}
&\nab^l Z_s^{t,\xi} = \nab^l \xi + \int_0^s \ind_{[t,T]}\nab^l Y_r^{t,\xi} \, dr \\
&\nab^l Y_s^{t,\xi} = \nab h\bigl(Z_T^{t,\xi}(\fdot)\bigr)\nab^l Z_T^{t,\xi} + 
\int_s^T \ind_{[t,T]}\nab F\bigl(r,Z_r^{t,\xi}(\fdot)\bigr)\nab^l Z_r^{t,\xi} \, dr 
 \notag \\
 & +\sum_{j=2}^l \nab^j h\bigl(Z_T^{t,\xi}(\fdot)\bigr)\Bigl[\sum_{i_1+\dots+i_j = l-j+1}\nab^{i_1}Z_T^{t,\xi}\ldots \nab^{i_j}Z_T^{t,\xi}\Bigr]
 \notag \\
 & +\int_s^T\ind_{[t,T]} \sum_{j=2}^{l} \nab^{j} F\bigl(r, Z_r^{t,\xi}(\fdot)\bigr)\Bigl[\sum_{i_1+\dots+i_j = l-j+1}\nab^{i_1}Z_r^{t,\xi}\ldots \nab^{i_j}Z_r^{t,\xi}\Bigr]\, dr
 \notag \\
 \label{bwd-dl-ext}
&  -\sqrt{2\nu} \int_s^T \nab^l X_r^{t,\xi} \,dW_r
 }
and note that its solution  
$(\nab^l Z^{t,\xi}_s, \nab^l Y^{t,\xi}_s, \nab^l X^{t,\xi}_s)$
can be obtained from the solution to
\eqref{fwd-dl}, \eqref{bwd-dl} by extending it to $[0,t]$
as follows: $\nab^l Z^{t,\xi}_s= \nab^l\xi$, 
$\nab^l Y^{t,\xi}_s = \nab^l Y^{t,\xi}_t$ , $\nab^l X^{t,\xi}_s = 0$, $s\in [0,t]$.
The extended triple solves the FBSDEs \eqref{fwd-dl-ext}, \eqref{bwd-dl-ext}
on $[0,T]$.
The same argument as in the proof of Lemma 14 of \cite{ns1}
implies that there exists a constant $\gm>0$ such that
\aa{
\|\nab^l Y^{t,\xi}_t - \nab^l Y^{t',\xi}_{t'}\|_{L_2(\Tor,\Rnu^{n^l})} \lt \gm |t-t'|.
}
Therefore  the map $t\mto Y^{t,\xi}_t$ is continuous with respect to
the $H^\al(\Tor,\Rnu^n)$-topology.
\end{proof}
\begin{lem}
\label{lem7}
Let the function $y: [0,T] \x\Tor \to \Rnu^n$ be defined by the formula
\aaa{
\label{function_y}
y(t,\te) = Y^{t,e}_t(\te).
}
Then, for every $t\in [0,T]$, there exists a set $\Om'\sub \Om$
of full $\PP$-measure, so that for all $u\in[t,T]$, for all $\om\in\Om'$
the following relation holds:
\aaa{
\label{111}
Y^{t,e}_u = y(u,\fdot) \circ Z^{t,e}_u.
}
\end{lem}
\begin{proof}
Note that \eqref{form_of_solution} implies that if $\xi$ is $\mc F_t$-measurable
then
\aaa{
\label{NS-rinv}
Y^{t,\xi}_t = y(t,\fdot)\circ \xi.
}
Further, for each fixed $u\in [t,T]$,  $(Z_s^{t,e}, Y_s^{t,e}, X_s^{t,e})$
 is a solution of the following problem on $[u,T]$:
 \[
 \begin{cases}
  Z_s^{t,e} = Z_u^{t,e} + \int_u^s Y^{t,e}_r dr + \sqrt{2\nu}\, (W_s - W_u),\\
  Y_s^{t,e} = h\bigl(Z_T^{t,e}(\fdot)\bigr) +  \int_s^T F\bigl(r,Z_r^{t,e}(\fdot)\bigr) dr - \sqrt{2\nu}\,\int_s^T X^{t,e}_r dW_r.
 \end{cases}
 \]
  By the uniqueness of solution, it holds that 
 $Y_s^{t,e} = Y_s^{u,Z_u^{t,e}}$ a.s. on $[u,T]$. Next, by \eqref{NS-rinv},
 we obtain that $Y_u^{u,Z_u^{t,e}}= y(u,\fdot) \circ Z_u^{t,e}$.
 This implies that there exists a set
 $\Om_u$ (which depends on $u$) of full $\PP$-measure
 such that \eqref{111} holds everywhere on $\Om_u$.
 Clearly, one can find a set $\Om_\Qnu$, $\PP(\Om_\Qnu)=1$, such that
 \eqref{111} holds on $\Om_\Qnu$ for all rational $u\in [t,T]$.
 But the trajectories of $Z^{t,e}_s$ and $Y^{t,e}_s$
 are a.s. continuous with  respect to $L_2(\Tor,\Rnu^n)$-topology
 by Lemma \ref{lem6}.
Furthermore $y(t,\fdot)$ is continuous in $t$ with respect to (at least) the 
$L_2(\Tor,\Rnu^n)$-topology. 
Therefore, \eqref{111} holds a.s. with respect to the $L_2(\Tor,\Rnu^n)$-topology. Since
 both sides of \eqref{111} are continuous in $\te\in\Tor$ it also holds a.s. for all $\te\in\Tor$.
\end{proof}
Finally the function 
$y(s,\fdot)$ defined by \eqref{function_y}
indeed verifies the Burgers equation. This is the content of the next lemma.
\begin{lem}
\label{y-smooth}

The function $y$ defined by formula \eqref{function_y}
is $C^1$-smooth in $t\in [0,T]$, and is a solution of problem \eqref{brgs}.
\end{lem}
\begin{proof}
Let $\dl > 0$. We obtain:
\aa{
y(t+\dl,\fdot) - y(t,\fdot) = Y^{t+\dl,e}_{t+\dl} - Y^{t,e}_t = Y^{t+\dl,e}_{t+\dl} -Y^{t,e}_{t+\dl} +  
Y^{t,e}_{t+\dl} - Y^{t,e}_t.
}
As before, let  $G^\al$ be the group of $H^\al$-diffeomorphisms $\Tor\to\Tor$,
and let $\hat Y_s$ be the right-invariant vector field
on $G^\al$ generated by $y(s,\fdot)$ (see \cite{ns1}). 
Relation \eqref{111} allows us to represent the SDE \eqref{fwd-e} as an SDE on the manifold
$G^\al$. Indeed, 
by results of \cite{Gliklikh} and \cite{ns1}, 
the SDE 
\aaa{
\label{sde-manifold}
\begin{cases}
dZ^{t,e}_s = \exp\{ \hat Y_s(Z^{t,e}_s)\, ds + \sqrt{2\nu}\, dW_s\},\\
Z^{t,e}_t = e
\end{cases}
}
where $\exp$ is the exponential map of the weak 
Riemannian metric on $G^\al$ (see \cite{ns1}),
has a unique $G^\al$-valued solution.
As it was proved in \cite{ns1}, 
the latter solution coincides with the unique solution of the $H^\al(\Tor,\Rnu^n)$-valued SDE
\aa{
\begin{cases}
dZ^{t,e}_s = y(s,\fdot)\circ Z^{t,e}_s\, ds + \sqrt{2\nu}\, dW_s,\\
Z^{t,e}_t = e.
\end{cases}
}
Therefore, the $Z^{t,e}_s$-part of the solution to
\eqref{fwd-e}, \eqref{bwd-e} 
is the unique solution to \eqref{sde-manifold}.
By Lemma \ref{lem7}, a.s.
$Y^{t,e}_{t+\dl} = \hat Y_{t+\dl}(Z^{t,e}_{t+\dl})$.
Thus we obtain that a.s.
\aa{
y(t+\dl,\fdot) - y(t,\fdot) = \bigl(\hat Y_{t+\dl}(e) - \hat Y_{t+\dl}(Z^{t,e}_{t+\dl})\bigr)
+(Y^{t,e}_{t+\dl} - Y^{t,e}_t).
}
We use the BSDE \eqref{bwd-e} for the second difference and
apply It\^o's formula to the first difference when considering
$\hat Y_{t+\dl}$ as a $C^2$-smooth function $G^\al \to L_2(\Tor,\Rnu^2)$. 
We obtain:
\mm{
\hat Y_{t+\dl}(Z^{t,e}_{t+\dl}) - \hat Y_{t+\dl}(e) 
 = \int_t^{t+\dl}  \hat Y^{t,e}_r(Z^{t,e}_r) [\hat Y_{t+\dl}(Z^{t,e}_r)] \, dr\\
+ \sqrt{2\nu} \, \int_t^{t+\dl} \sum_{i=1}^n
\bigl[\bar\nab_{e_i} \,\hat Y_{t+\dl}(Z^{t,e}_r)\bigr] \, dW_r
+ 2\nu \, \int_t^{t+\dl} \sum_{i=1}^n\bigl[\bar\nab^2_{e_i} \hat Y_{t+\dl}(Z^{t,e}_r)\bigr] \, dr
}
where $\bar\nab$ is the covariant derivative on $G^\al$,
$e_i$ are regarded as constant vector fields on $G^\al$,
and the expression
$\hat Y^{t,e}_r(Z^{t,e}_r) [\hat Y_{t+\dl}(Z^{t,e}_r)]$ has the same meaning as 
in \cite{ns1}.
We obtain:
\mm{
\hat Y_{t+\dl}(Z^{t,e}_{t+\dl}) - \hat Y_{t+\dl}(e) =  \int_t^{t+\dl}\hspace{-2mm} dr \, \nab_{y(r,\fdot)}\,y(t+\dl,\fdot)\circ Z^{t,e}_r \\
+\int_t^{t+\dl} \hspace{-2mm} dr \, \nu\,\lap\, y(t+\dl,\fdot)\circ Z^{t,e}_r 
+ \sqrt{2\nu} \,\int_t^{t+\dl} 
\sum_{i=1}^n
\bigl[\bar\nab_{e_i} \,\hat Y_{t+\dl}(Z^{t,e}_r)\bigr] \, dW_r.
}
Further we have:
\aa{
Y^{t,e}_{t}- Y^{t,e}_{t+\dl} = \int_t^{t+\dl}\hspace{-2mm} dr \, F(r,\fdot) \circ
Z^{t,e}_r - \sqrt{2\nu}\,\int_t^{t+\dl} \hspace{-2mm} X^{t,e}_r \, dW_r.
}
Taking expectations implies:
\mmm{
\label{diff}
\frac1{\dl}\, \bigl(y(t+\dl,\fdot) - y(t,\fdot)\bigr) = - \frac1{\dl}\, \E \Bigl[
\int_t^{t+\dl}\hspace{-2mm} dr \, [ \, (y(r,\fdot),\nab)\,y(t+\dl,\fdot)\\
+\nu\,\lap\, y(t+\dl,\fdot)+F(r,\fdot)] \circ Z^{t,e}_r
\Bigr].
}
Note that $Z^{t,e}_r$, $F(r,\fdot)$, and 
$(y(r,\fdot),\nab)\,y(t+\dl,\fdot)\circ Z^{t,e}_r$
are continuous in $r$ a.s. with respect to the $L_2(\Tor,\Rnu^2)$-topology.
By Lemma \ref{lem4},
$\nab\, y(t,\fdot)$ and $\lap\, y(t,\fdot)$ are continuous in $t$ with
respect to at least the $L_2(\Tor,\Rnu^2)$-topology. 
Formula \eqref{diff} and the fact that $Z^{t,e}_t = e$ imply that in the $L_2(\Tor,\Rnu^2)$-topology
\aaa{
\label{ns3}
\pl_t y(t,\fdot) = -[\nab_{y(t,\fdot)}\,y(t,\fdot)
+\nu\,\lap\, y(t,\fdot) + F(t,\fdot)]. 
}
Since the right-hand side of \eqref{ns3} is an $H^{\al-2}$-map, 
so is the left-hand side. This implies that $\pl_t y(t,\fdot)$ is continuous in $\te\in \Tor$.
Therefore, \eqref{ns3} holds for any $\te\in\Tor$.
Relation \eqref{ns3} is obtained so far for the right derivative of $y(t,\te)$
with respect to $t$.
Note that the right-hand side of \eqref{ns3} is continuous in $t$ which implies that
the right derivative $\pl_t y(t,\te)$ is continuous in $t$ on
$[0,T)$. Hence, it is uniformly continuous on every compact subinterval of $[0,T)$.
This implies the existence of the left derivative of $y(t,\te)$ in $t$, and therefore,
the existence of the continuous derivative $\pl_t y(t,\te)$ everywhere on $[0,T]$.
\end{proof}
\begin{rem}
{\rm
Note that at the same time we have proved that the process $Z^{t,e}_s$ takes values in 
the group $G^\al$ of $H^\al$-diffeomorphisms $\Tor\to\Tor$.
} 
\end{rem}
\section*{Acknowledgements} 
The first author acknowledges the support of the Portuguese Foundation for Science and Technology through the project PTDC/MAT/69635/2006.
The second author acknowledges the support of the Portuguese Foundation for Science and Technology through  
the Centro de Matem\'atica da Universidade do Porto.


\begin{thebibliography}{99}
\bibitem[A]{}V.I.~Arnold, \it Sur la g\'eom\'etrie
 diff\'erentielle des groupes de Lie de dimension infinie et ses applications
a l'hydrodynamique des fluides parfaits, \rm Ann. Inst. Fourier 16 (1966), 316--361.
 %
 \bibitem[B]{belopolskaya}
        {Ya. I. Belopolskaya, Yu. L. Dalecky},  \textit{Stochastic equations
        and differential geometry}, Series: Mathematics and its Applications,
        Kluwer Academic Publishers, Netherlands, (1989), 260 p.
 %
\bibitem[C-S-T-V]{}P.~ Cheridito, H. ~Mete Soner, N.~Touzi  and N.~Victoir, 
\it Second order backward stochastic differential equations and fully nonlinear parabolic PDEs, \rm  Comm. Pure Appl. Math.   60  (2007),  1081--1110.       
        
%
\bibitem[C-S]{ns1}A.B. Cruzeiro, E. Shamarova, \textit{Navier--Stokes equations and forward–backward 
SDEs on the group of diffeomorphisms of a torus},
Stochastic Processes and their Applications, 119, (2009), 4034--4060
%
\bibitem[D]{delarue}F. Delarue, \textit{On the existence and uniqueness of solutions to
the FBSDEs in a non-generate case}, Stoch. Proc. and their Appl.
\textbf{99}, (2002), 209--286.
%
\bibitem[G]{gliklikh4}
 Yu.\ E.\ Gliklikh,
\textit{Solutions of Burgers, Reynolds, and Navier--Stokes equations 
via stochastic perturbations of inviscid flows}, 
Journal of Nonlinear Mathematical Physics,
Vol. 17, No. Supplementary Issue 1 (2010) 15–-29.
 %
 \bibitem[G1]{Gliklikh}
 Yu.\ E.\ Gliklikh,
 \textit{Global Analysis in Mathematical Physics: Geometric and Stochastic Methods},
 Springer (1997), 213 p.
%
\bibitem[E-M]{}D.G.~Ebin and J.~Marsden, 
\it Groups of diffeomorphisms and the motion of an incompressible fluid, \rm  Ann. of Math.   92  (1970),  102--163.
%
\bibitem[N-Y-Z]{}T.~Nakagomi, K.~Yasue, J.-C.~Zambrini, \it
 Stochastic variational derivations of the Navier-Stokes equation, \rm Lett. Math. Phys., 160  (1981), 337--365.
%
\bibitem[Y]{}K.~Yasue, \it
 A variational principle for the Navier-Stokes equation,\rm J. Funct. Anal., 51 (2), (1983),
  133--141
\end{thebibliography}
\end{document}